\documentclass[12pt,twoside]{article}
\usepackage[english]{babel}
\usepackage[latin1]{inputenc}
\usepackage{amsmath}
\usepackage{amssymb,amsfonts}
\usepackage{graphicx}                   

\newcommand{\bb}{\mathbf{b}}

\newcommand{\bt}{\mathbf{t}}

\newcommand{\HYP}{\mathbb{H}^3}
\newcommand{\HYN}{\mathbb{H}^n}

\begin{document}
\pagestyle{myheadings}
\markboth{\centerline{Jen\H o Szirmai}}
{Density upper bound for congruent and non-congruent hyperball packings $\ldots$}
\title
{Density upper bound for congruent and non-congruent hyperball packings generated by truncated regular simplex tilings}

\author{\normalsize{Jen\H o Szirmai} \\
\normalsize Budapest University of Technology and \\
\normalsize Economics Institute of Mathematics, \\
\normalsize Department of Geometry \\
\date{\normalsize{\today}}}

\maketitle


\begin{abstract}
In this paper we study congruent and non-congruent hyperball (hypersphere) packings of the truncated regular tetrahedron tilings.
These are derived from the Coxeter simplex tilings $\{p,3,3\}$ $(7\le p \in \mathbb{N})$ and $\{5,3,3,3,3\}$
in $3$ and $5$-dimensional hyperbolic space.
We determine the densest hyperball packing arrangements related to the above tilings. We find packing densities
using congruent hyperballs and determine the smallest density upper bound of non-congruent hyperball packings generated by the above tilings.
\end{abstract}

\newtheorem{theorem}{Theorem}[section]
\newtheorem{corollary}[theorem]{Corollary}
\newtheorem{conjecture}{Conjecture}[section]
\newtheorem{lemma}[theorem]{Lemma}
\newtheorem{exmple}[theorem]{Example}
\newtheorem{defn}[theorem]{Definition}
\newtheorem{rmrk}[theorem]{Remark}
\newenvironment{definition}{\begin{defn}\normalfont}{\end{defn}}
\newenvironment{remark}{\begin{rmrk}\normalfont}{\end{rmrk}}
\newenvironment{example}{\begin{exmple}\normalfont}{\end{exmple}}
\newenvironment{acknowledgement}{Acknowledgement}


\section{Introduction}
The classical sphere packing problems concern arrangements of non-overlap\-ping equal spheres which fill a space,
usually the three dimensional Euclidean space. However, ball (sphere) packing problems can be generalized to the other
$3$-dimensional Thurston geometries and to higher various dimensional spaces, but the difficulty lies -- similar to the case of hyperbolic space --
in the rigorous definition of the packing density.

In $n$-dimensional space of constant curvature $\mathbb{E}^n$, $\mathbb{H}^n$, $\mathbb{S}^n$ $(n \ge 2)$ let
$d_n(r)$ be the density of $n+1$ spheres
of radius $r$ mutually touching one another with respect to the simplex spanned by the centres of the spheres.
L.~Fejes T\'oth and H.~S.~M.~Coxeter
conjectured that in an $n$-dimensional space of constant curvature the packing density balls of radius $r$ can not exceed $d_n(r)$.
This conjecture has been proved by C.~Rogers in the Euclidean space $\mathbb{E}^n$ \cite{Ro64}.
The 2-dimensional spherical case was settled by L.~Fejes T\'oth in \cite{FTL}. %
In the $n$-dimensional hyperbolic space $\mathbb{H}^n$ $(n\ge2)$ there are $3$ kinds
of ''balls (spheres)": balls (spheres), horoballs (horospheres) and hyperballs (hyperspheres).

K.~B\"or\"oczky in \cite{B78} proved that the above conjecture holds for balls (spheres) and for horoballs if the horoballs of the
{\it same type} (for the notion of horoball type see \cite{Sz12}, \cite{Sz12-2}).
\begin{rmrk}
\item[1.] In hyperbolic space $\mathbb{H}^3$ this result can be extended to $r=\infty$ \cite{B78} where the densest
horoball packing can be realized by different regular arrangements \cite{KSz}.
\item[2.] If we allow horoballs of different types at the various vertices of a totally asymptotic simplex and
generalize the notion of the simplicial density function in $\mathbb{H}^n$, $(n \ge 2)$ then the B\"or\"oczky--Florian
type density upper bound does not remain valid for the fully asymptotic simplices \cite{Sz12}, \cite{Sz12-2}.
\item[3.] In the papers \cite{KSz14} and \cite{Sz15} we investigated the horoball packings related to the $4$-dimensional Coxeter simplices and
to the hyperbolic $24$-cell. The densest horoball packing arrangements of the a above tilings provide $\approx 0.71645$
density which is the densest known ball and horoball packing density in $\mathbb{H}^4$.
\end{rmrk}

In \cite{Sz06-1} and \cite{Sz06-2} we studied the regular prism tilings and the corresponding optimal hyperball packings in
$\mathbb{H}^n$ $(n=3,4)$ and in the paper  \cite{Sz13-3} we extended the method developed in previous papers 
to 5-dimensional hyperbolic space, and constructed to each investigated Coxeter tiling a corresponding regular prism tiling,
and studied their optimal hyperball packings by congruent hyperballs,
moreover, we determined their metric data and densities. \newline

In the hyperbolic plane $\mathbb{H}^2$ the universal upper bound of the hypercycle packing density is $\frac{3}{\pi}$
as proved by I.~Vermes in \cite{V79}. Recently, (to the author's best knowledge) the candidates for the densest hyperball
(hypersphere) packings in the $3,4$ and $5$-dimensional hyperbolic space $\mathbb{H}^n$ are derived by the regular prism
tilings studied in papers \cite{Sz06-1}, \cite{Sz06-2} and \cite{Sz13-3}.

In $\mathbb{H}^2$ the universal lower bound of the hypercycle covering density is $\frac{\sqrt{12}}{\pi}$
determined again by I.~Vermes in \cite{V81}.
In the paper \cite{Sz13-4} we studied the $n$-dimensional $(n \ge 3)$ hyperbolic regular prism honeycombs
and the corresponding coverings by congruent hyperballs and we determined their least dense covering densities.
Moreover, we have formulated a conjecture for the candidate of the least dense hyperball
covering by congruent hyperballs in the 3- and 5-dimensional hyperbolic space.

In \cite{Sz14} we studied the problem of hyperball (hypersphere) packings in
$3$-dimensional hyperbolic space.
We described a procedure for each saturated hyperball packing to get a decomposition of 3-dimensional hyperbolic space $\HYP$ into truncated
tetrahedra. Therefore, in order to get a density upper bound for hyperball packings it is sufficient to determine
the density upper bound of hyperball packings in truncated simplices.
Thus, we considered the hyperball packings in truncated simplices and proved that if the truncated tetrahedron is regular, then the
density of the densest packing is $\approx 0.86338$ which is larger than the B\"oröczky-Florian density upper bound,
however these hyperball packing configurations are only locally optimal and cannot be extended to the entirety
of the hyperbolic spaces $\mathbb{H}^3$.
Moreover, we proved that the known densest hyperball packing related to regular prism tilings can be realized by a regular truncated tetrahedral tiling
\cite{Sz06-1}.

In this paper we study the problem of congruent and non-congruent hyperball (hypersphere) packings to each truncated regular tetrahedron tiling.
These are derived from the Coxeter simplex tilings $\{p,3,3\}$ and $\{5,3,3,3,3\}$ in the $3$ and $5$-dimensional hyperbolic space.
We determine the densest hyperball packing arrangements and their densities
with congruent hyperballs and determine the smallest density upper bound of non-congruent hyperball packings generated by the above tilings.
\section{Basic notions}
\subsection{Complete orthoschemes}
For $\mathbb{H}^n$ we use the projective model in the Lorentz space $\mathbb{E}^{1,n}$ of signature $(1,n)$,
i.e.~$\mathbb{E}^{1,n}$ denotes the real vector space $\mathbf{V}^{n+1}$ equipped with the bilinear
form of signature $(1,n)$, that is
$
\langle ~ \mathbf{x},~\mathbf{y} \rangle = -x^0y^0+x^1y^1+ \dots + x^n y^n
$
where the non-zero vectors
$
\mathbf{x}=(x^0,x^1,\dots,x^n)\in\mathbf{V}^{n+1} \ \  \text{and} \ \ \mathbf{y}=(y^0,y^1,\dots,y^n)\in\mathbf{V}^{n+1},
$
are determined up to real factors, for representing points of $\mathcal{P}^n(\mathbb{R})$. Then $\mathbb{H}^n$ can be interpreted
as the interior of the quadric
$
Q=\{[\mathbf{x}]\in\mathcal{P}^n | \langle ~ \mathbf{x},~\mathbf{x} \rangle =0 \}=:\partial \mathbb{H}^n
$
in the real projective space $\mathcal{P}^n(\mathbf{V}^{n+1},
\mbox{\boldmath$V$}\!_{n+1})$.

The points of the boundary $\partial \mathbb{H}^n $ in $\mathcal{P}^n$
are called points at infinity of $\mathbb{H}^n $, the points lying outside $\partial \mathbb{H}^n $
are said to be outer points of $\mathbb{H}^n $ relative to $Q$. Let $P([\mathbf{x}]) \in \mathcal{P}^n$, a point
$[\mathbf{y}] \in \mathcal{P}^n$ is said to be conjugate to $[\mathbf{x}]$ relative to $Q$ if
$\langle ~ \mathbf{x},~\mathbf{y} \rangle =0$ holds. The set of all points which are conjugate to $P([\mathbf{x}])$
form a projective (polar) hyperplane
$
pol(P):=\{[\mathbf{y}]\in\mathcal{P}^n | \langle ~ \mathbf{x},~\mathbf{y} \rangle =0 \}.
$
Thus the quadric $Q$ induces a bijection
(linear polarity $\mathbf{V}^{n+1} \rightarrow
\mbox{\boldmath$V$}\!_{n+1})$)
from the points of $\mathcal{P}^n$
onto its hyperplanes.

The point $X [\bold{x}]$ and the hyperplane $\alpha [\mbox{\boldmath$a$}]$
are called incident if $\bold{x}\mbox{\boldmath$a$}=0$ ($\bold{x} \in \bold{V}^{n+1} \setminus \{\mathbf{0}\}, \ \mbox{\boldmath$a$} \in \mbox{\boldmath$V$}_{n+1}
\setminus \{\mbox{\boldmath$0$}\}$).
\begin{definition}
An orthoscheme $\mathcal{S}$ in $\mathbb{H}^n$ $(2\le n \in \mathbb{N})$ is a simplex bounded by $n+1$ hyperplanes $H^0,\dots,H^n$
such that
(see \cite{B--H, K89})
$
H^i \bot H^j, \  \text{for} \ j\ne i-1,i,i+1.
$
\end{definition}

{\it The orthoschemes of degree} $m$ in $\mathbb{H}^n$ are bounded by $n+m+1$ hyperplanes
$H^0,H^1,\dots,H^{n+m}$ such that $H^i \perp H^j$ for $j \ne i-1,~i,~i+1$, where, for $m=2$,
indices are taken modulo $n+3$. For a usual orthoscheme we denote the $(n+1)$-hyperface opposite to the vertex $A_i$
by $H^i$ $(0 \le i \le n)$. An orthoscheme $\mathcal{S}$ has $n$ dihedral angles which
are not right angles. Let $\alpha^{ij}$ denote the dihedral angle of $\mathcal{S}$
between the faces $H^i$ and $H^j$. Then we have
$
\alpha^{ij}=\frac{\pi}{2}, \ \ \text{if} \ \ 0 \le i < j -1 \le n.
$
The $n$ remaining dihedral angles $\alpha^{i,i+1}, \ (0 \le i \le n-1)$ are called the
essential angles of $\mathcal{S}$.
Geometrically, complete orthoschemes of degree $d$ can be described as follows:
\begin{enumerate}
\item
For $m=0$, they coincide with the class of classical orthoschemes introduced by
{{Schl\"afli}} (see Definitions 2.1).
The initial and final vertices, $A_0$ and $A_n$ of the orthogonal edge-path
$A_iA_{i+1},~ i=0,\dots,n-1$, are called principal vertices of the orthoscheme.
\item
A complete orthoscheme of degree $m=1$ can be interpreted as an
orthoscheme with one outer principal vertex, say $A_n$, which is truncated by
its polar plane $pol(A_n)$ (see Fig.~1 and 3). In this case the orthoscheme is called simply truncated with
outer vertex $A_n$.
\item
A complete orthoscheme of degree $m=2$ can be interpreted as an
orthoscheme with two outer principal vertices, $A_0,~A_n$, which is truncated by
its polar hyperplanes $pol(A_0)$ and $pol(A_n)$. In this case the orthoscheme is called doubly
truncated. We distinguish two different types of orthoschemes but 
will not enter into the details (see \cite{K89}).
\end{enumerate}
An $n$-dimensional tiling $\mathcal{P}$ (or solid tessellation, honeycomb) is an infinite set of
congruent polyhedra (polytopes) that fit together to fill all space $(\mathbb{H}^n~ (n \geqq 3))$ exactly once,
so that every face of each polyhedron (polytope) belongs to another polyhedron as well.
At present the cells are congruent orthoschemes and later regular truncated tetrahedra.
A tiling with orthoschemes exists if and only if each dihedral angle of a tile is submultiple of $2\pi$
(in the hyperbolic plane the zero angle is also possible).

Another approach to describing tilings involves the analysis of their symmetry groups.
If $\mathcal{P}$ is such a simplex tiling, then any motion taking one cell into another maps the
entire tiling onto itself. The symmetry group of this tiling is denoted by
$Sym \mathcal{P}$.
Therefore the simplex is a fundamental domain of the group $Sym \mathcal{P}$ generated by reflections in its
$(n-1)$-dimensional hyperfaces.

For the schemes (weighted graph) of {\it complete Coxeter orthoschemes} $\mathcal{S} \subset \mathbb{H}^n$
we adopt the usual conventions and sometimes even use them in the Coxeter case:
Two nodes are joined by an edge if the corresponding hyperplanes are not orthogonal.
If two nodes are related by the weight $\cos{\frac{\pi}{p}}$
then they are joined by a ($p-2$)-fold line for $p=3,~4$ and by a single line marked $p$ for $p \geq 5$.
In the hyperbolic case if two bounding hyperplanes of $S$ are parallel, then the corresponding nodes
are joined by a line marked $\infty$. If they are divergent then their nodes are joined by a dotted line.

The ordered set $\{k_1,\dots,k_{n-1},k_n\}$ is said to be the
Coxeter-Schl$\ddot{a}$fli symbol of the simplex tiling $\mathcal{P}$ generated by $\mathcal{S}$.
To every scheme there is a corresponding
symmetric matrix $(c^{ij})$ of size $(n+1)\times(n+1)$ where $c^{ii}=1$ and, for $i \ne j\in \{0,1,2,\dots,n \}$,
$c^{ij}$ equals $-\cos{\frac{\pi}{k_{ij}}}$ with all angles between the facets $i$,$j$ of $\mathcal{S}$.

For example, $(c^{ij})$ below is the so called Coxeter-Schl\"afli matrix of the orthoscheme $S$ in
3-dimensional hyperbolic space $\mathbb{H}^3$ with
parameters (nodes) $k_1=p,k_2=q,k_3=r$ :
\[
(c^{ij}):=\begin{pmatrix}
1& -\cos{\frac{\pi}{p}}& 0 & 0 \\
-\cos{\frac{\pi}{p}} & 1 & -\cos{\frac{\pi}{q}}& 0 \\
0 & -\cos{\frac{\pi}{q}} & 1 & -\cos{\frac{\pi}{r}} \\
0 & 0 & -\cos{\frac{\pi}{r}} & 1 \\
\end{pmatrix}. \tag{2.1}
\]
In general the complete Coxeter orthoschemes were classified by {{Im Hof}} in
\cite{IH85} by generalizing the method of {{Coxeter}} and {{B\"ohm}}, who
showed that they exist only for dimensions $\leq 9$. From this classification it follows, that the complete
orthoschemes of degree $m=1$ exist up to 5 dimensions.

In this paper we consider some tilings generated by orthoschemes of degree $1$ where the initial vertex $A_n$ is outer point regarding the quadric $Q$. 
These orthoschemes and the corresponding Coxeter tilings exist in the $2$-, $3-$, $4-$ and
$5-$dimensional hyperbolic spaces and
are characterized by their Coxeter-Schl\"afli symbols and graphs.

In $n$-dimensional hyperbolic space $\mathbb{H}^n$ $(n \ge 2)$
it can be seen that if $\mathcal{O}=A_0A_1A_2 \dots A_n$ $P_0P_1P_2 \dots P_n$ is a complete
orthoscheme  with degree $d=1$ (a simply frustum orthoscheme) where $A_n$ is a outer vertex of
$\mathbb{H}^n$ then the points $P_0,P_1,P_2,\dots,P_{n-1}$ lie on the polar hyperplane $\pi$ of $A_n$ (see Fig.~1 in $\HYP$).
\begin{figure}[ht]
\centering
\includegraphics[width=6.5cm]{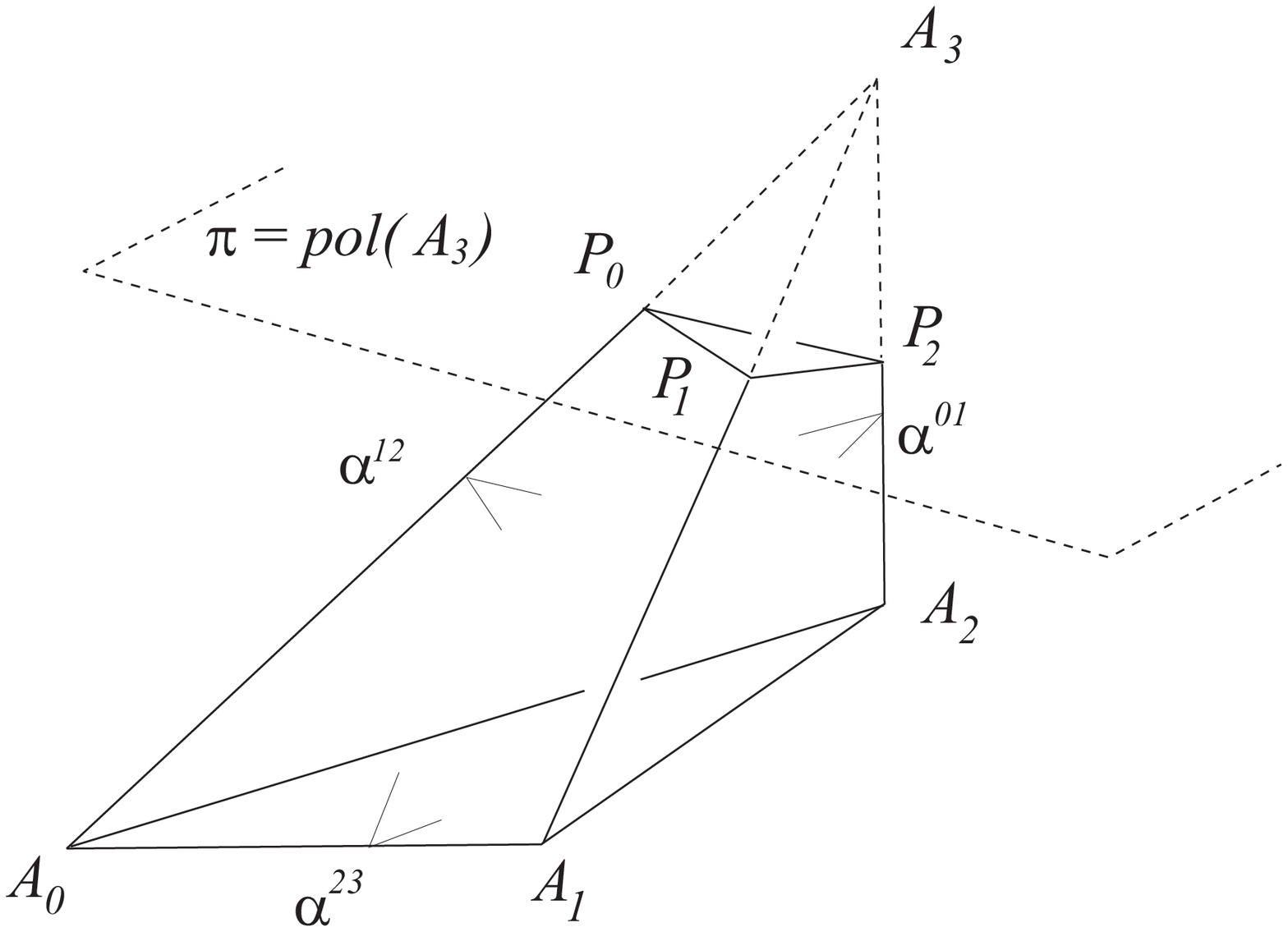} \includegraphics[width=6.5cm]{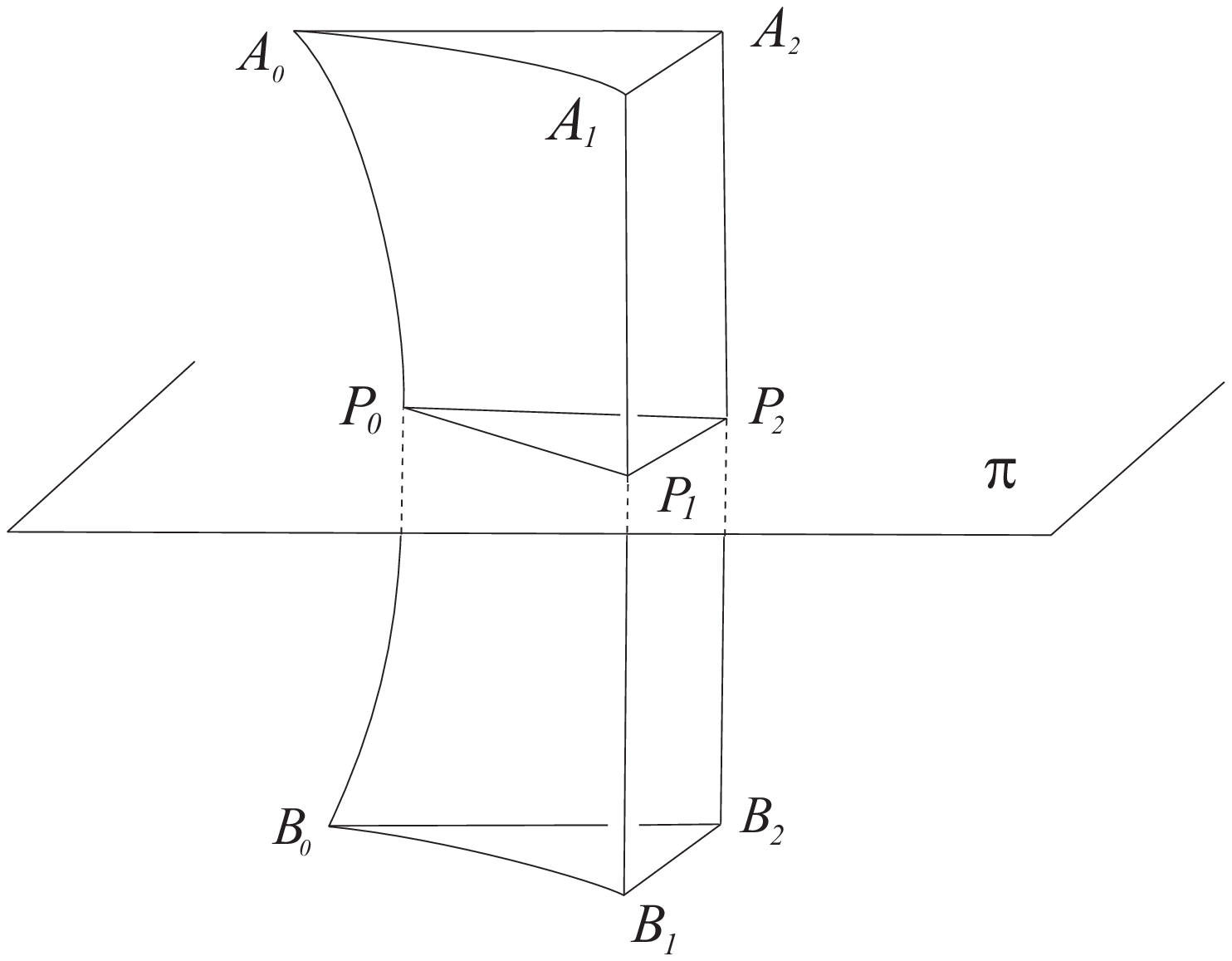}

a. \hspace{5cm} b.
\caption{a.~A $3$-dimensional complete orthoscheme of degree $m=1$ (simple frustum orthoscheme)
with outer vertex $A_3$. This orthoscheme is truncated by its polar plane $\pi=pol(A_3)$.~b.~Two congruent
adjacent simple frustum orthoschemes.}
\label{}
\end{figure}
The images of $\mathcal{O}$ under reflections on its side facets generate a tiling in $\HYN$.

{\it The constant $k =\sqrt{\frac{-1}{K}}$ is the natural length unit in
$\mathbb{H}^n$. $K$ will be the constant negative sectional curvature. In the following we assume that $k=1$.}
\subsubsection{Volumes of the $n$-dimensional \\ Coxeter orthoschemes}
The truncated regular simplex tilings exist only in $3-$ and $5-$dimensional hyperbolic spaces (derived from the Coxeter tilings $\{p,3,3\}$
and $\{5,3,3,3,3\}$) therefore
we use now the $3-$ and $5-$dimensional volume formulas of the above truncated Coxeter orthoschemes.
\begin{enumerate}
\item $3$-dimensional hyperbolic space $\HYP$:
{Our polyhedron $A_0A_1A_2P_0P_1P_2$ is a simple frustum orthoscheme with
outer vertex $A_3$ (see Fig.~1) whose volume can be calculated by the following theorem of R.~Kellerhals
\cite{K89}:}
\begin{theorem} The volume of a three-dimensional hyperbolic
complete ortho\-scheme (except Lambert cube cases) $\mathcal{S}$
is expressed with the essential angles $\alpha_{01},\alpha_{12},\alpha_{23}, \ (0 \le \alpha_{ij} \le \frac{\pi}{2})$
(Fig.~1) in the following form:
\begin{align}
&Vol_3(\mathcal{O})=\frac{1}{4} \{ \mathcal{L}(\alpha_{01}+\theta)-
\mathcal{L}(\alpha_{01}-\theta)+\mathcal{L}(\frac{\pi}{2}+\alpha_{12}-\theta)+ \notag \\
&+\mathcal{L}(\frac{\pi}{2}-\alpha_{12}-\theta)+\mathcal{L}(\alpha_{23}+\theta)-
\mathcal{L}(\alpha_{23}-\theta)+2\mathcal{L}(\frac{\pi}{2}-\theta) \}, \tag{2.2}
\end{align}
where $\theta \in [0,\frac{\pi}{2})$ is defined by the following formula:
$$
\tan(\theta)=\frac{\sqrt{ \cos^2{\alpha_{12}}-\sin^2{\alpha_{01}} \sin^2{\alpha_{23}
}}} {\cos{\alpha_{01}}\cos{\alpha_{23}}}
$$
and where $\mathcal{L}(x):=-\int\limits_0^x \log \vert {2\sin{t}} \vert dt$ \ denotes the
Lobachevsky function.
\end{theorem}
For our prism tilings $\mathcal{T}_{pqr}$ we have:~
$\alpha_{01}=\frac{\pi}{p}, \ \ \alpha_{12}=\frac{\pi}{q}, \ \
\alpha_{23}=\frac{\pi}{r}$ .
\item $5$-dimensional hyperbolic space $\mathbb{H}^5$:

{{R.~Kellerhals}} in \cite{K92} developed a procedure to determine the volumes of 5-dimensional hyperbolic orthoschemes,
moreover the volumes of the complete orthoschemes $\mathcal{S}_{pqrst}$ with Coxeter-Schl\"afli symbol $[5,3,3,3,3]$ and
$[5,3,3,3,4]$ can be computed by the differential volume formula of L.~ Schl\"afli:
\begin{equation*}
Vol_5(\mathcal{O}_{pqrst})=\frac{1}{4}\int_{\alpha_i}^{\frac{2\pi}{5}}{Vol_3([5,3,\beta(t)])dt+\frac{\zeta(3)}{3200}}
\end{equation*}
with a compact tetrahedron $[5,3,\beta(t)]$ whose angle parameter $0<\beta(t)<\frac{\pi}{2}$ is given by
$$
\beta({t})=\arctan{\sqrt{2-\cot^2{t}}}.
$$
Then, the volume of the 3-dimensional orthoscheme face $[5,3,\beta(t)]$ as given by Lobachevsky's formula:
\begin{equation}
\begin{gathered}
Vol_3([5,3,\beta(t)])=\frac{1}{4} \{\mathcal{L}_2\big(\frac{\pi}{5}+\theta(t)\big)-\mathcal{L}_2\big(\frac{\pi}{5}-\theta(t)\big)-\mathcal{L}_2\big(\frac{\pi}{6}+\theta(t)\big)+\\
\mathcal{L}_2\big(\frac{\pi}{6}-\theta(t)\big)+\mathcal{L}_2\big(\beta(t)+\theta(t)\big)-\mathcal{L}_2\big(\beta(t)-\theta(t)\big)+2\mathcal{L}_2\big(\frac{\pi}{2}-\theta(t)\big) \end{gathered} \tag{2.3}
\end{equation}
where $\mathcal{L}(\omega)$ is the Lobachevsky's function,
$
\theta(t)=\arctan\frac{\sqrt{1-4\sin^2\frac{\pi}{5} \sin^2{\beta(t)}}}{2\cos\frac{\pi}{5} \cos\beta(t)}
$ and $\beta(t)=\frac{\pi}{3}$ or $\frac{\pi}{4}$.
\end{enumerate}
\subsection{On hyperballs}
The equidistant surface (or hypersphere) is a quadratic surface that lies at a constant distance
from a plane in both halfspaces. The infinite body of the hypersphere is called hyperball.
The $n$-dimensional {\it half-hypersphere } $(n=3,5)$ with distance $h$ to a hyperplane $\pi$
is denoted by $\mathcal{H}_n^{h+}$.
The volume of a bounded hyperball piece $\mathcal{H}_n^{h+}(\mathcal{A}_{n-1})$
bounded by an $(n-1)$-polytope $\mathcal{A}_{n-1} \subset \pi$, $\mathcal{H}_n^{h+}$ and by
hyperplanes orthogonal to $\pi$ derived from the facets of $\mathcal{A}_{n-1}$ can be determined by the
formulas (2.4) and (2.5) that follow from the suitable extension of the classical method of {{J.~Bolyai}}:
\begin{equation}
Vol_3(\mathcal{H}_3^{h+}(\mathcal{A}_2))=\frac{1}{4}Vol_2(\mathcal{A}_{2})\left[\sinh{2h}+
2 h \right], \tag{2.4}
\end{equation}
\begin{equation}
Vol_5(\mathcal{H}_5^{h+}(\mathcal{A}_4))=\frac{1}{16} Vol_4(\mathcal{A}_4)\left[\left( \frac{1}{2} \sinh{4h}+
4 \sinh{2h}\right) +6 h \right], \tag{2.5}
\end{equation}
where the volume of the hyperbolic $(n-1)$-polytope $\mathcal{A}_{n-1}$ lying in the plane $\pi$ is $Vol_{n-1}(\mathcal{A}_{n-1})$.
\section{On hyperball packings in truncated regular simplices}
In \cite{Sz14} we studied the problem of congruent hyperball (hypersphere) packings in 
$3$-dimensional hyperbolic space. We described to each saturated congruent hyperball packing a procedure to get a decomposition
of $3$-dimensional hyperbolic space $\HYP$ into truncated
tetrahedra. Therefore, in order to get a density upper bound to hyperball packings it is sufficient to determine
the density upper bound of hyperball packings in truncated simplices.

Moreover, we proved that the {\it known densest congruent hyperball packingn} relating to the regular prism tilings can be
realized by a regular truncated tetrahedron tiling \cite{Sz06-1}, \cite{Sz14}.

Similarly to the above question it is interesting to construct and to study the locally optimal congruent and {\it non-congruent} hyperball packings
relating to truncated regular simplex tilings in $3$- and higher dimensions as well.

{\it Therefore, the $n$-dimensional regular simplex tilings and the corresponding optimal congruent and non-congruent hyperball packings
(similarly to Euclidean space and to hyperbolic horosphere packings) play an important role among the packing problems.}

In $n$-dimensional ($2 < n \in \mathbb{N}$) hyperbolic spaces there are two types of the regular simplex tilings which are derived
from the Coxeter simplex tilings $\{p,3,3\}$ $(\mathbb{N} \ni p \ge 7)$ and $\{5,3,3,3,3\}$ in $3$ and $5$-dimensional hyperbolic space.

\subsection{Hyperball packings with congruent hyperballs in $3$-dimensional regular truncated tetrahedra}
We consider an arbitrary saturated congruent hyperball packing $\mathcal{B}^{h(p)}$ of hyperballs $\mathcal{H}^{h_i(p)}$. In \cite{Sz14} we showed
that there is a decomposition of $\HYP$ into truncated tetrahedra.
One of them $\mathcal{S}(p)$=$C_1^1 C_2^1 C_3^1$ $C_1^2 C_2^2 C_3^2$ $C_1^3 C_2^3 C_3^3$ $C_1^4 C_2^4 C_3^4$
is illustrated in Fig.~2-3.

The ultraparallel base planes of $\mathcal{H}^{h_i(p)}$ $(i=1,2,3,4)$ are denoted by $\beta_i$. Their poles $B_i$ are outer points in the Cayley-Klein model
(see Fig.~2). The distance between two base
planes $d(\beta_i,\beta_j)=:e_{ij}$ ($i < j$, $i,j \in \{1,2,3,4\})$ and $d$ is the hyperbolic distance function)
at least $2h(p)$. Moreover, the volume of the truncated simplex $\mathcal{S}(p)$ is denoted by $Vol(\mathcal{S}(p))$.
We introduce the locally density function $\delta(\mathcal{S}(h(p)))$ related to $\mathcal{S}(p)$:
\begin{definition}
\begin{equation}
\delta(\mathcal{S}(h(p))):=\frac{\sum_{i=1}^4 Vol(\mathcal{H}^{h_i(p)} \cap \mathcal{S}(p))}{Vol({\mathcal{S}(p)})}. \tag{3.1}
\end{equation}
\end{definition}
In the following we assume that the ultraparallel base planes $\beta_i$ of $\mathcal{H}^{h(p)}_i$ $(i=1,2,3,4)$
generate a "regular truncated tetrahedron" $\mathcal{S}^r$ with outer vertices $B_i$ (see Fig.~3.~a) i.e.
the non-orthogonal dihedral angles of $\mathcal{S}^r$ are equal to $\frac{2\pi}{p}$,
$(7 \le p\in \mathbb{N})$ and the distance between two base planes $d(\beta_i,\beta_j)=:e_{ij}$  ($i < j$, $i,j \in \{1,2,3,4\})$ are equal to $2h(p)$.
In this case for a given parameter $p$ the length of the common perpendiculars $h(p)=\frac{1}{2}e_{ij}$ $(i < j$, $i,j \in \{1,2,3,4\})$
can be determined by the machinery of the projective geometry.

The points $P_2[{\mathbf{p}}_2]$ and $Q_2[{\mathbf{q}}_2]$ are proper points of hyperbolic $3$-space and
$Q_2$ lies on the polar hyperplane $pol(B_1)[\mbox{\boldmath$b$}^1]$ of the outer point $B_1$ thus
\begin{equation}
\begin{gathered}
\mathbf{q}_2 \sim c \cdot \mathbf{b}_1 + \mathbf{p}_2 \in \mbox{\boldmath$b$}^1 \Leftrightarrow
c \cdot \mathbf{b}_1 \mbox{\boldmath$b$}^1+\mathbf{p}_2 \mbox{\boldmath$b$}^1=0 \Leftrightarrow
c=-\frac{\mathbf{p}_2 \mbox{\boldmath$b$}^1}{\mathbf{b}_1 \mbox{\boldmath$b$}^1} \Leftrightarrow \\
\mathbf{q}_2 \sim -\frac{\mathbf{p}_2 \mbox{\boldmath$b$}^1}{\mathbf{b}_1 \mbox{\boldmath$b$}^1}
\mathbf{b}_1+\mathbf{p}_2 \sim \mathbf{p}_2 (\mathbf{b}_1 \mbox{\boldmath$b$}^1) - \mathbf{b}_1 (\mathbf{p}_2 \mbox{\boldmath$b$}^1)=
\mathbf{p}_2 h_{33}-\mathbf{b}_1 h_{23},
\end{gathered} \tag{3.2}
\end{equation}
where $h_{ij}$ is the inverse of the Coxeter-Schl\"afli matrix ($i,j=0,1,2,3)$)
\[(c^{ij}):=\begin{pmatrix}
1& -\cos{\frac{\pi}{p}}& 0&0\\
-\cos{\frac{\pi}{p}} & 1 & -\cos{\frac{\pi}{3}}&0\\
0 & -\cos{\frac{\pi}{3}} & 1&-\cos{\frac{\pi}{3}} \\
0&0&-\cos{\frac{\pi}{3}}&1\\
\end{pmatrix} \notag
\]
of the orthoscheme $\mathcal{O}$.
The hyperbolic distance $h(p)$ can be calculated by the following formula:
\[
\begin{gathered}
\cosh{h(p)}=\cosh{P_2Q_2}=\frac{- \langle {\mathbf{q}}_2, {\mathbf{p}}_2 \rangle }
{\sqrt{\langle {\mathbf{q}}_2, {\mathbf{q}}_2 \rangle \langle {\mathbf{p}}_2, {\mathbf{p}}_2 \rangle}}= \\ =\frac{h_{23}^2-h_{22}h_{33}}
{\sqrt{h_{22}\langle \mathbf{q}_2, \mathbf{q}_2 \rangle}} =
\sqrt{\frac{h_{22}~h_{33}-h_{23}^2}
{h_{22}~h_{33}}}.
\end{gathered} \tag{3.3}
\]
We get that the volume $Vol(\mathcal{S}^r)$, the maximal height $h(p)$ of the congruent hyperballs lying in $\mathcal{S}^r$ and the
$\sum_{i=1}^4 Vol(\mathcal{H}^h_i \cap \mathcal{S}^r))$ depend only on the parameter $p$ of the truncated regular tetrahedron $\mathcal{S}^r$.

Therefore, the density $\delta(\mathcal{S}^r(h(p)))$ is depended only on parameter $p$ $(7 \ge p\in \mathbb{N})$. Moreover,
the volume of the hyperball pieces can be computed by the formula (2.4) and the volume of $\mathcal{S}^r$ can be determined by the Theorem 2.2.

In \cite{Sz14} we proved the following
\begin{theorem}
The density function $\delta(\mathcal{S}^r(h(p)))$, $p\in (6,\infty)$
is attained its maximum at $p^{opt} \approx 6.13499$ and $\delta(\mathcal{S}^r(h(p)))$
is strictly increasing on the interval $(6,p^{opt})$ and strictly decreasing on the interval $(p^{opt},\infty)$. Moreover, the optimal density
$\delta^{opt}(\mathcal{S}^r(h(p^{opt}))) \approx 0.86338$ which is larger than the B\"oröczky-Florian density upper bound,
however these hyperball packing configurations are only locally optimal and cannot be extended to the entirety
of the hyperbolic spaces $\mathbb{H}^3$.
\end{theorem}
\begin{figure}[ht]
\centering
\includegraphics[width=8cm]{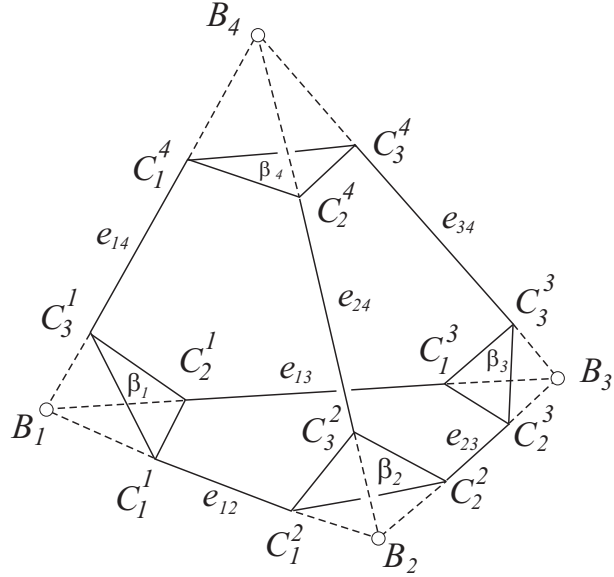}
\caption{Truncated tetrahedron}
\label{}
\end{figure}
\begin{corollary}
 The density function $\delta(\mathcal{S}^r(h(p)))$, $\mathbb{N} \ni p \ge 7$ attains its maximum at the parameter
$p=7$. Congruent hyperball packings $\mathcal{B}^{h(p)}$ with the above parameters related to the regular truncated tetrahedra
can be extended to the entire hyperbolic space. The maximal density is $\delta(\mathcal{S}^r(h(7))) \approx 0.82251$.
\end{corollary}
\begin{rmrk}We note here that these coincide with the hyperball packings to the regular prism tilings in $\HYP$ with Schl\"afli symbols
$\{p,3,3\}$ which are discussed in \cite{Sz06-1} because their vertex figure is tetrahedron given by Schl\"afli symbol \{3,3\}.
\end{rmrk}
\begin{figure}[ht]
\centering
\includegraphics[width=5.5cm]{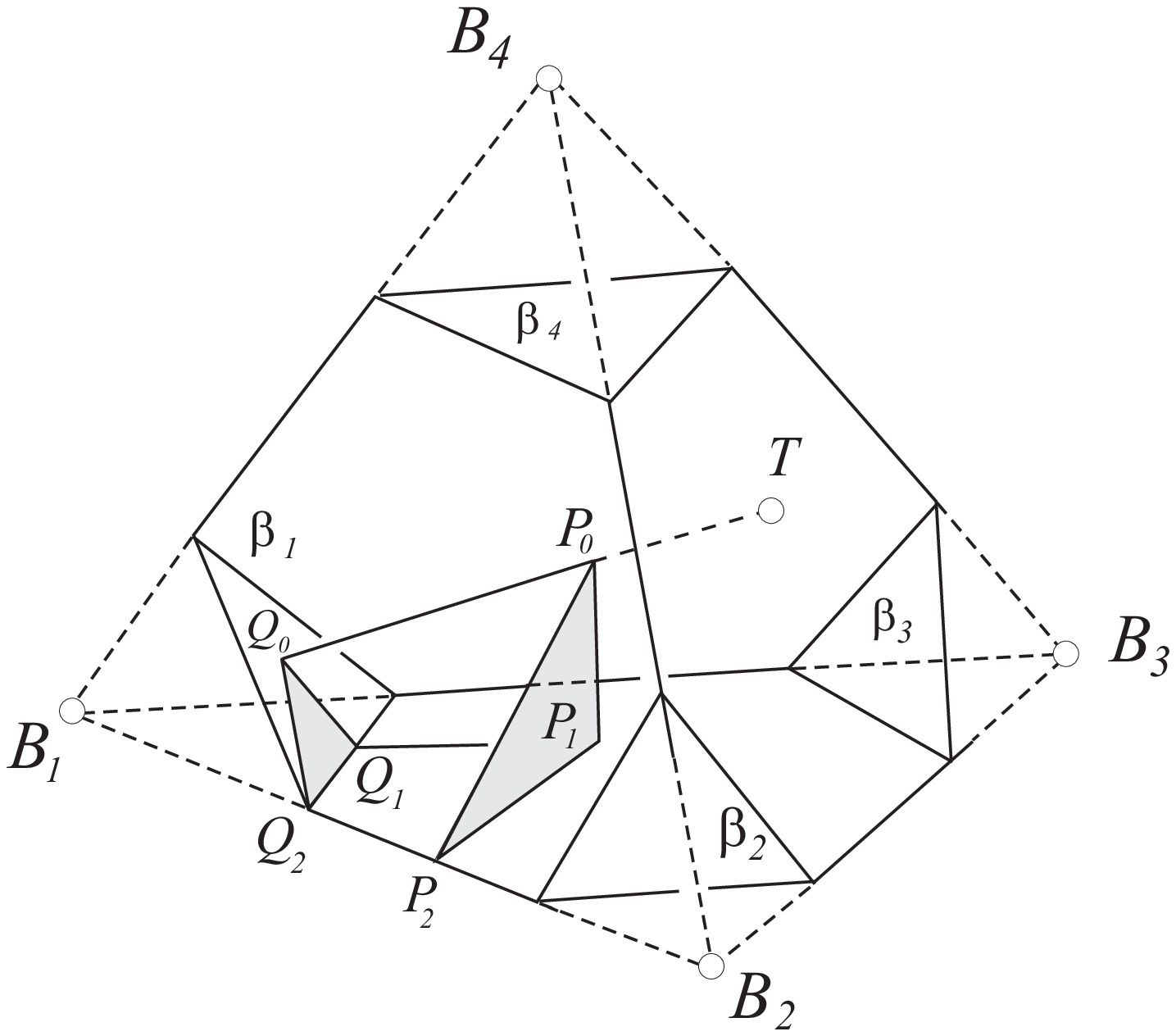} \includegraphics[width=6.5cm]{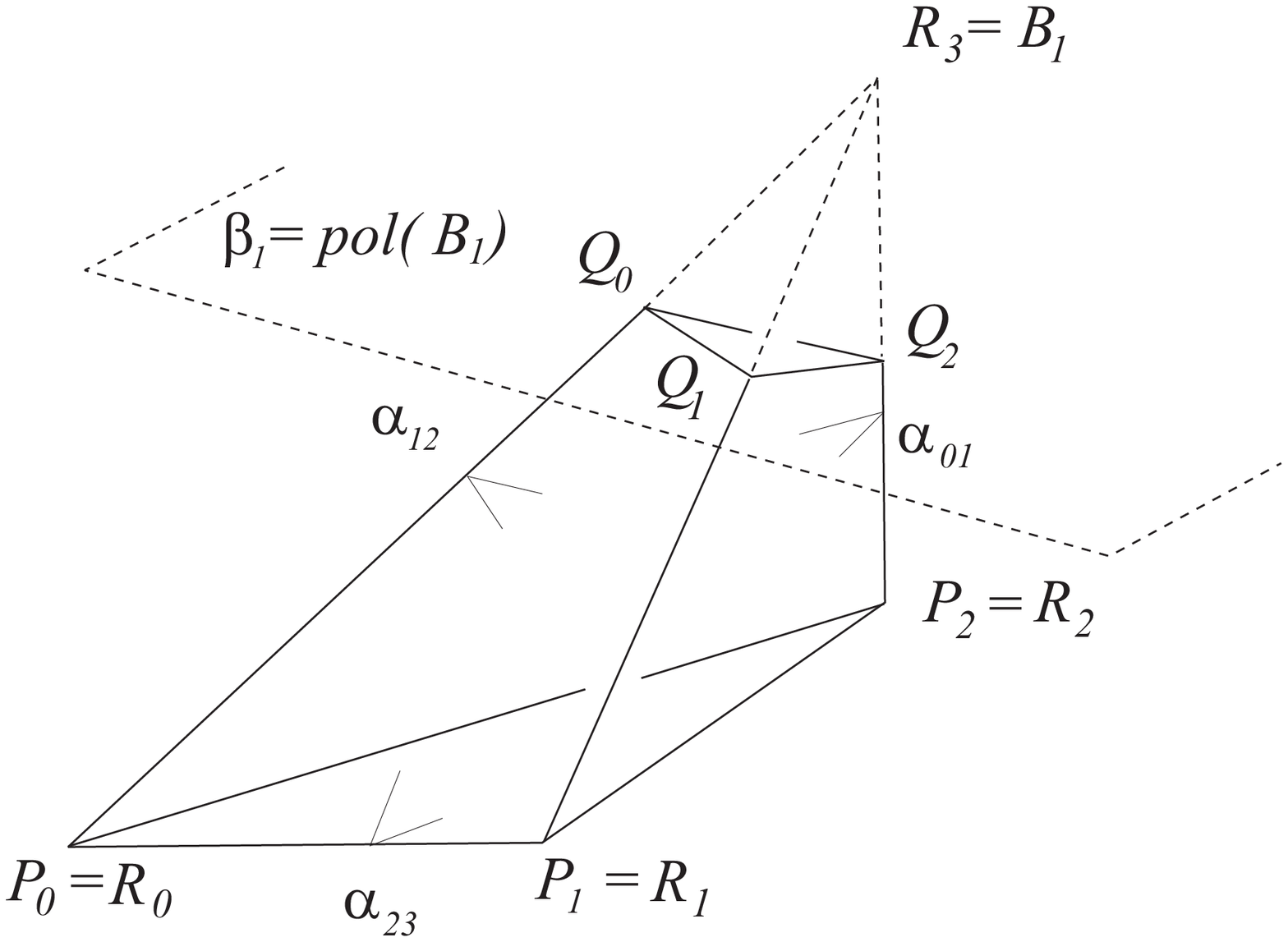}

a.~~~~~~~~~~~~~~~~~~~~~~~~~~~~~~~~~~~~~b.
\caption{Truncated tetrahedron with a complete orthoscheme of degree $m=1$ (simple frustum orthoscheme)}
\label{}
\end{figure}
In the following Table we summarize the data of the hyperball packings for some parameters $p$, ($7 \ge p\in \mathbb{N}$).
\medbreak
\centerline{\vbox{
\halign{\strut\vrule~\hfil $#$ \hfil~\vrule
&\quad \hfil $#$ \hfil~\vrule
&\quad \hfil $#$ \hfil\quad\vrule
&\quad \hfil $#$ \hfil\quad\vrule
&\quad \hfil $#$ \hfil\quad\vrule
\cr
\noalign{\hrule}
\multispan5{\strut\vrule\hfill\bf Table 1, \hfill\vrule}%
\cr
\noalign{\hrule}
\noalign{\vskip2pt}
\noalign{\hrule}
p & h(p) & Vol(\mathcal{O}) & Vol(\mathcal{H}^{h(p)}(\mathcal{A}))& \delta^{opt} \cr
\noalign{\hrule}
7 & 0.78871 & 0.08856 & 0.07284 & 0.82251 \cr
\noalign{\hrule}
8 & 0.56419 & 0.10721 & 0.08220 & 0.76673 \cr
\noalign{\hrule}
9 & 0.45320 & 0.11825 & 0.08474 & 0.71663 \cr
\noalign{\hrule}
\vdots & \vdots  & \vdots  & \vdots  & \vdots \cr
\noalign{\hrule}
20 & 0.16397 & 0.14636 & 0.06064 & 0.41431 \cr
\noalign{\hrule}
\vdots & \vdots  & \vdots  & \vdots  & \vdots \cr
\noalign{\hrule}
50 & 0.06325 & 0.15167 & 0.02918 & 0.19240 \cr
\noalign{\hrule}
\vdots & \vdots  & \vdots  & \vdots  & \vdots \cr
\noalign{\hrule}
100 & 0.03147 & 0.15241 & 0.01549 & 0.10165 \cr
\noalign{\hrule}
p \to \infty & 0 & 0.15266 & 0 & 0 \cr
\noalign{\hrule}}}}
\subsection{Hyperball packings with non-congruent hyperballs in $3$-dimensional regular truncated tetrahedra}
We consider a regular truncated tetrahedron tiling $\mathcal{T}(\mathcal{S}^r(p))$ with Schl\"afli symbol $\{p,3,3\}$, $(7 \le p\in\mathbb{N})$.
One tile of it $\mathcal{S}^r(p)$=$C_1^1 C_2^1 C_3^1$ $C_1^2 C_2^2 C_3^2$ $C_1^3 C_2^3 C_3^3$ $C_1^4 C_2^4 C_3^4$
is illustrated in Fig.~2. The ultraparallel base planes of the {\it not absolutely congruent hyperballs} (hyperspheres) $\mathcal{H}^{h_i}_i$
$(i=1,2,3,4)$ are denoted by $\beta_i$ e.g.
the base plane $\beta_i$ of the hyperball $\mathcal{H}^{h_i}_i$ (with height $h_i$) is determined by vertices $C_1^iC_2^iC_3^i$
(see Fig.~2, $i=1,2,3,4$). The distance between two base
planes $d(\beta_i,\beta_j)=:e_{ij}$ are equal ($i < j$, $i,j \in \{1,2,3,4\})$.
Moreover, the volume of the truncated simplex $\mathcal{S}^r(p)$ is denoted by $Vol(\mathcal{S}^r(p))$, similarly to the above section.

The centre of the hexagonal side face opposite the vertex $B_i$ of $\mathcal{S}^r(p)$, (it is rectangular hexagon)
is denoted by $T_i$ ($i\in \{1,2,3,4\}$) (see Fig.~2 and Fig.~3.a).
We construct non-congruent hyperball packings to $\mathcal{T}(\mathcal{S}^r(p))$ tilings therefore the hyperballs have to satisfy the
following requirements:
\begin{enumerate}
\item The base plane $\beta_i$ of the hyperball $\mathcal{H}^{h_i}_i$ (with height $h_i$) is determined by vertices $C_1^iC_2^iC_3^i$ (see Fig.~2),
\item $ card \{ int (\mathcal{H}^{h_i}_i)\cap int(\mathcal{H}^{h_j}_j)\}=0$, $i\ne j$,
\item $e_{ij} \le h_i+h_j$, $i,j\in \{1,2,3,4\}, ~ i<j$,
\item $ card \{ int (\mathcal{H}^{h_i}_i \cap int(B_jB_kB_l ~\text{plane})\}=0$ $(i,j,k,l \in \{1,2,3,4\},~ i~\ne j,k,l,~ j<k<l$ i.e.
the distance $w(p):=d(\beta_i,T_i)\ge h_i$.
\end{enumerate}
{\it If the hyperballs hold the above requirements then we obtain hyperball packings $\mathcal{B}(\mathcal{S}^r(p))$ in hyperbolic $3$-space
derived by the structure of the considered Coxeter simplex tilings.}

We introduce the locally density function $\delta(\mathcal{S}^r(p))$ related to $\mathcal{S}^r(p)$:
\begin{definition}
\begin{equation}
\delta(\mathcal{S}^r(p)):=\frac{\sum_{i=1}^4 Vol(\mathcal{H}^{h_i}_i \cap \mathcal{S}^r(p))}{Vol({\mathcal{S}^r(p)})}. \tag{3.4}
\end{equation}
\end{definition}
It is well known that a packing is locally optimal (i.e. its density is locally maximal), then it is locally stable i.e. each ball is fixed by
the other ones so that no ball of packing
can be moved alone without overlapping another ball of the given ball packing or by other requirements of the corresponding tiling.

We set up from the optimal congruent ball arrangement $\mathcal{B}^{h(p)}$ (see former section)
where the congruent hyperballs touch each other at the
"midpoints" of the edges of $\mathcal{S}^r(p)$. We choose an arbitrary hyperball (e.g. $\mathcal{H}^{h(p)}_1$) and blow up
this hyperball (hypersphere) keeping the hyperballs $\mathcal{H}^{h_i}_i$ $(i=2,3,4)$ tangent to it upto this hypersphere
\begin{enumerate}
\item touches the faces $B_2B_3B_4$ at point $T_1$, here the height of the hyperball $\mathcal{H}^{w(p)}_1$ is $w(p)$ if $w(p) \le 2h(p)$
(it touches its opposite face). During this expansion the height of $\mathcal{H}^{h_1}_1$ is $h_1=h(p)+x$ $(x \in \mathbb{R}^+)$ where
$h(p)\le h(p)+ x \le w(p)$. The height of further hyperballs are $h_2=h_3=h_4=h(p)-x$. (If $x=0$ then the horoballs are congruent.)
\item passing through the point $C_1^2$ (touches the planes $\beta_i$ $i=2,3,4$), here the height of the extended
hyperball $\mathcal{H}^{2h(p)}_1$ is $2h(p)$ if $2h(p) \le w(p)$ (it touches the planes $\beta_i$).
During this expansion the height of $\mathcal{H}^{h_1}_1$ is $h_1=h(p)+x$ $(x \in \mathbb{R}^+)$ where
$h(p)\le h(p)+ x \le 2h(p)$, the height of further hyperballs are $h_2=h_3=h_4=h(p)-x$. (If $x=0$ then the horoballs are congruent.)
\end{enumerate}
We extend this procedure to images of the hyperballs $\mathcal{H}^{h_i}_i$ $(i=1,2,3,4)$ by the considered Coxeter group and obtain
non-congruent hyperball arrangements $\mathcal{B}^x(p)$.

{\it The main problem is: what is the maximum of density function $\delta(\mathcal{S}^r(x,p))$ for a given integer parameter $p \ge 7$ where
$x \in \mathbb{R}$, and $x \in [0,h(p)]$ or $x \in [0,w(p)-h(p)]$}.

During this expansion process we can compute for a given integer parameter $p \ge 7$ the densities
$\delta(\mathcal{S}^r(x,p))$ of considered packings as the function of $x$.
\subsubsection{Computations for parameter $p \ge 7$}
Every $n$-dimensional hyperbolic truncated regular simplex can be derived from a $n$-dimensional regular Euclidean simplex.
We introduce a projective coordinate system (see Section 2) and a unit sphere $\mathbb{S}^{n-1}$ centred at the origin which is
interpreted as the ideal boundary of $\overline{\mathbb{H}}^n$ in Beltrami-Cayley-Klein's ball model.

Now, we consider a $3$-dimensional regular Euclidean tetrahedron centred at the origin with outer vertices regarding the Beltrami-Cayley-Klein's ball model.
The projective coordinates of the vertices of this tetrahedron are
\begin{equation}
\begin{gathered}
B_1\Big(1,\frac{2\sqrt{2}y}{3},0,\frac{-y}{3}\Big); ~ B_2\Big(1,\frac{-\sqrt{2}y}{3},\frac{\sqrt{2}y}{\sqrt{3}},\frac{-y}{3}\Big);~
B_3\Big(1,\frac{-\sqrt{2}y}{3},\frac{\sqrt{2}y}{\sqrt{3}},\frac{-y}{3}\Big);\\
B_{4}(1,0,0,y); ~ ~\text{where} ~ ~ 1 < y\in \mathbb{R}.
\end{gathered} \tag{3.5}
\end{equation}
The truncated tetrahedron $\mathcal{S}^r(p)$ can be derived from the above tetrahedron by cuttings with the polar planes of 
vertices $B_i$ $(i=1,2,3,4)$. The images of $\mathcal{S}^r(p)$ under reflections on its side facets generate a tiling in $\HYP$
if its non-right dihedral angles are $\frac{2\pi}{p}$ $(\mathbb{N} \ni p\ge 7)$. It is easy to see, that if the parameter $p$ is given, then
\begin{equation}
y=\sqrt{3}\sqrt{\frac{3\cos{\frac{2\pi}{p}}-1}{\cos{\frac{2\pi}{p}}+1}}. \tag{3.6}
\end{equation}
We have to determine for any parameter $p$ the distances $h(p)$ and $w(p)$. The values of $h(p)$ can be derived
from formula (3.3) and $w(p)$  from the next formula:
\begin{equation}
\begin{gathered}
\sinh{w(p)}=\Bigg|{\frac{\langle \bb_1,\bt_1\rangle}{\sqrt{-\langle \bb_1,\bb_1\rangle \langle \bt_1,\bt_1\rangle}}}\Bigg|=
\frac{y^2+3}{3\sqrt{(y^2-1)(-y^2+9)}}. \tag{3.7}
\end{gathered}
\end{equation}
If $p=7$ then we obtain the following results:
\begin{equation}
\begin{gathered}
2h(7)\approx 1.57741; w(7) \approx 1.51843 \Rightarrow ~w(7)<2h(7) ~ \Rightarrow \\
\Rightarrow ~x \in [0, w(7)-h(7)\approx 0.72972].
\end{gathered} \tag{3.8}
\end{equation}
We note here, that if $x=0$ then the hyperspheres are congruent (see the former section).
Therefore, we can compute during the expansion process for the given integer parameter $p =7$ the densities of
$\delta(\mathcal{S}^r(x,7))$ (see Definition 3.4) of considered packings as the function of $x$ using the formulas (2.4), (3.6), (3.7), (3.8) and Theorem 2.2:
\begin{equation}
\begin{gathered}
\delta(\mathcal{S}^r(x,7))=\frac{Vol(\mathcal{H}^{h(7)+x} \cap \mathcal{S}^r(x,7))+3 \cdot Vol(\mathcal{H}^{h(7)-x}
\cap \mathcal{S}^r(x,7))}{Vol({\mathcal{S}^r(x,7)})}, \\
\text{where} ~ x \in [0, w(7)-h(7)\approx 0.72972].
\end{gathered} \tag{3.9}
\end{equation}
The graph of $\delta(\mathcal{S}^r(x,7))$ is described in Fig.~4. Analyzing the above density function we get that the
maximal density is achieved at the starting point of the above interval (at the congruent case)
with density $\approx 0.82251$ (see Table 1). The density in the endpoint of the above interval is
$\delta(\mathcal{S}^r(w(7)-h(7),7)) \approx 0.74649$.
\begin{figure}[ht]
\centering
\includegraphics[width=8cm]{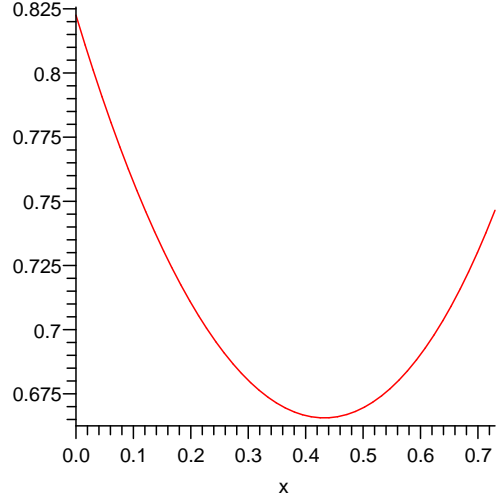}
\caption{Density function of $\delta(\mathcal{S}^r(x,7))$ where $x \in [0, w(7)-h(7)\approx 0.72972]$.}
\label{}
\end{figure}
Similarly to the above computations for parameter $p=7$ we can analyze the density functions and their maximums of
non-congruent hyperball packings generated by considered truncated simplex tilings (or Coxeter tilings) for all possible integer parameters
$p \ge 7$. Using the results of Theorem 3.2 and Corollary 3.3 we obtain the following
\begin{theorem}
\begin{enumerate}
\item The maximum of the density function $\delta(\mathcal{S}^r(x,p))$ is attained at the starting
point of the interval $x \in [0,h(p)]$ or $x \in [0,w(p)-h(p)]$ depending on integer parameter $p\ge 7$ i.e.
the congruent hyperball packing provides the densest hyperball packing for a given parameter $p$.
\item The maximum of the density function $\delta(\mathcal{S}^r(x,p))$ $(p \ge 7$ integer parameter) is achieved at the parameters $x=0$, $p=7$.
Therefore, the density upper bound of the congruent and non-congruent hyperball packings is $\approx 0.82251$.
\end{enumerate}
\end{theorem}
\begin{rmrk}
We note here that the density function $\delta(\mathcal{S}^r(x,p))$, $p\in (6,\infty)$
is attained its maximum at $p^{opt} \approx 6.13499$ similarly to the congruent cases with density $\approx 0.86338$.

Therefore there are infinitely many non-congruent hyperball packing whose density
larger than the B\"oröczky-Florian density upper bound,
however these hyperball packing configurations are only locally optimal and cannot be extended to the entirety
of the hyperbolic spaces $\mathbb{H}^3$.
\end{rmrk}
\subsection{Hyperball packings with congruent and non-congruent hyperballs related to the $5$-dimensional regular truncated simplex tiling}
We consider the regular truncated 5-dimensional simplex tiling $\mathcal{T}(\mathcal{S}^r)$ generated by the Coxeter group with
Schl\"afli symbol $\{5,3,3,3,3\}$.
The vertices of a tile from the above tiling are given
in a projective coordinate system (see Section 2) where a unit sphere $\mathbb{S}^{4}$ centred at the origin provides
the ideal boundary of $\overline{\mathbb{H}}^5$ in Beltrami-Cayley-Klein's ball model.

The above regular truncated tetrahedron can be derived from a regular Euclidean tetrahedron
with outer vertices, regarding the Beltrami-Cayley-Klein's ball model. The
projective coordinates of its vertices are the following:
\begin{equation}
\begin{gathered}
B_1(1,\frac{y}{\sqrt{15}},\frac{y}{\sqrt{10}},\frac{y}{\sqrt{6}},\frac{y}{\sqrt{3}},y); ~
B_2(1,\frac{y}{\sqrt{15}},\frac{y}{\sqrt{10}},\frac{y}{\sqrt{6}},\frac{y}{\sqrt{3}},-y),\\
B_3(1,\frac{y}{\sqrt{15}},\frac{y}{\sqrt{10}},\frac{-\sqrt{3}y}{\sqrt{2}},0,0),~
B_4(1,\frac{y}{\sqrt{15}},\frac{-2\sqrt{2}y}{\sqrt{5}},0,0,0),\\
B_5(1,\frac{y}{\sqrt{15}},\frac{y}{\sqrt{10}},\frac{y}{\sqrt{6}},\frac{-2y}{\sqrt{3}},0),
B_6(1,\frac{-\sqrt{5}y}{\sqrt{3}},0,0,0,0) ~ ~\text{where} ~ ~ \frac{\sqrt{3}}{\sqrt{5}} < y\in \mathbb{R}.
\end{gathered} \tag{3.10}
\end{equation}
The truncated tetrahedron $\mathcal{S}^r$ can be derived from the above simplex by cuttings with the polar hyperplanes of
vertices $B_i$ $(i=1,2,3,4,5,6)$. Its non-right dihedral angles are $\frac{2\pi}{5}$.
The images of $\mathcal{S}^r$ under reflections on its side facets generate a tiling in $\mathbb{H}^5$.

It is easy to compute the parameter $y \in \mathbb{R}$:
\begin{equation}
\cos{\frac{2\pi}{5}}=\frac{y^2+3}{-y^2+15} \Rightarrow y \approx 1.11769. \tag{3.11}
\end{equation}

Similarly to the $3$-dimensional cases the ultraparallel base hyperplanes of the {\it not absolutely congruent hyperballs} (hyperspheres)
$\mathcal{H}^{h_i}_i$ $(i=1,2,3,4,5,6)$ are denoted by $\beta_i$.
The distance between two base
planes $2 h:=d(\beta_i,\beta_j)$ are equal for any $i < j$, $i,j \in \{1,2,3,4,5,6\}$, ($d$ is the hyperbolic distance function).
Moreover, the volume of the truncated simplex $\mathcal{S}^r$ is denoted by $Vol(\mathcal{S}^r)$, similarly to the above section.

We introduce the locally density function $\delta(\mathcal{S}^r)$ related to $\mathcal{S}^r$:
\begin{definition}
\begin{equation}
\delta(\mathcal{S}^r):=\frac{\sum_{i=1}^6 Vol(\mathcal{H}^{h_i}_i \cap \mathcal{S}^r)}{Vol({\mathcal{S}^r})}. \notag
\end{equation}
\end{definition}
The centre of the facet opposite the vertex $B_i$ of $\mathcal{S}^r$, is denoted by $T_i$ ($i\in \{1,2,3,4,5,6\}$) and
the distance between the hyperplane $\beta_i$ and $T_i$ is denoted by $w:=d(\beta_i,T_i)$.

Similarly to the $3$-dimensional computations we determine for given parameters the distances $h$ and $w$.
\begin{enumerate}
\item The values of $h$ can be derived
from the generalization of the formula (3.3):
\[
\begin{gathered}
\cosh{h}=\sqrt{\frac{h_{55}~h_{66}-h_{56}^2}
{h_{55}~h_{66}}} \approx 0.38360 \Rightarrow d(\beta_i,\beta_j)=2h \approx 0.76720.
\end{gathered} \tag{3.12}
\]
where $h_{ij}$ is the inverse of the Coxeter-Schl\"afli matrix $(c^{ij})$ $(i,j=0,1,2,3,4,5)$ of the orthoscheme $\mathcal{O}$.
\[(c^{ij}):=\begin{pmatrix}
1& -\cos{\frac{\pi}{5}}& 0&0&0&0 \\
-\cos{\frac{\pi}{5}} & 1 & -\cos{\frac{\pi}{3}}&0&0&0\\
0 & -\cos{\frac{\pi}{3}} & 1&-\cos{\frac{\pi}{3}},0&0 \\
0&0&-\cos{\frac{\pi}{3}}&1&-\cos{\frac{\pi}{3}}&0\\
0&0&0&-\cos{\frac{\pi}{3}}&1&-\cos{\frac{\pi}{3}}\\
0&0&0&0&-\cos{\frac{\pi}{3}}&1
\end{pmatrix}. \notag
\]
\item $w$ can be determined by the following formula:
\begin{equation}
\begin{gathered}
\sinh{w}=\Bigg|{\frac{\langle \bb_i,\bt_i\rangle}{\sqrt{-\langle \bb_i,\bb_i\rangle \langle \bt_i,\bt_i\rangle}}}\Bigg|=
\frac{\sqrt{5}(y^2+3)}{3\sqrt{(5y^2-3)(-y^2+15)}} \Rightarrow \\ \Rightarrow w \approx 1.15080. \tag{3.13}
\end{gathered}
\end{equation}
\end{enumerate}
We construct congruent and non-congruent hyperball packings to $\mathcal{T}(\mathcal{S}^r)$ tiling
therefore the hyperballs have to satisfy the following requirements:
\begin{enumerate}
\item The base plane $\beta_i$ of the hyperball $\mathcal{H}^{h_i}_i$ (with height $h_i$) is determined by the polar hyperplane
of vertex $B_i$,
\item $ card \{ int (\mathcal{H}^{h_i}_i)\cap int(\mathcal{H}^{h_j}_j)\}=0$,  $i\ne j$,
\item $e_{ij} \ge h_i+h_j$, $i,j\in \{1,2,3,4,5,6\}, ~ i<j$,
\item $ card \{ int (\mathcal{H}^{h_i}_i) \cap int(B_jB_kB_lB_sB_t ~\text{plane})\}=0$, \\ 
$i,j,k,l,s,t \in \{1,2,3,4,5,6\}$,~ $i~\ne j,k,l,s,t,~ 
j<k<l<s<t$.
\end{enumerate}
{\it If the hyperballs hold the above requirements then hyperball packings $\mathcal{B}(\mathcal{S}^r)$ can be derived by the
corresponding Coxeter group in $5$-dimensional hyperbolic space.}

We set up from the optimal congruent ball arrangement $\mathcal{B}^h$ 
where the congruent hyperballs touch each other at the
"midpoints" of the edges of $\mathcal{S}^r$. We choose an arbitrary hyperball (e.g. $\mathcal{H}^{h}_1$) and blow up
this hyperball (hypersphere) keeping the hyperballs $\mathcal{H}^{h_i}_i$ $(i=2,3,4,5,6)$ tangent to it upto this hypersphere
touches the planes $\beta_i$, (see (3.12), (3.13)). Here the height of the extended
hyperball $\mathcal{H}^{2h}_1$ is $2h$ because $2h(p) < w$.
During this expansion the height of $\mathcal{H}^{h_1}_1$ is $h_1=h+x$ $(x \in \mathbb{R}^+)$ where
$h \approx 0.38360 \le h+ x \le 2h \approx 0.76720$, the height of further hyperballs are $h_2=h_3=h_4=h_5=h_6=h-x$.

We extend this procedure to images of the hyperballs $\mathcal{H}^{h_i}_i$ $(i=1,2,3,4,5,6)$ by the considered Coxeter group and obtain
non-congruent hyperball packings $\mathcal{B}^x$.

{\it The main problem is: what is the maximum of density function $\delta(\mathcal{S}^r(x))$ where
$x \in \mathbb{R}$, and $x \in [0,h]$}.

We note here, that if $x=0$ then the hyperspheres are congruent and if $x \in (0,h]$ then we obtain non-congruent hyperball packings.

Therefore, we can compute during the expansion process the densities of
$\delta(\mathcal{S}^r(x))$ (see Definition 3.7) of considered packings as the function of $x$
using the formulas (2.5), (3.11), (3.12), (3.13) and (2.3):
\begin{equation}
\begin{gathered}
\delta(\mathcal{S}^r(x))=\frac{Vol(\mathcal{H}^{h+x} \cap \mathcal{S}^r(x))+5 \cdot Vol(\mathcal{H}^{h-x}
\cap \mathcal{S}^r(x))}{Vol({\mathcal{S}^r(x)})}, \\
\text{where} ~ x \in [0, h \approx 0.38360].
\end{gathered} \tag{3.14}
\end{equation}
The graph of $\delta(\mathcal{S}^r(x))$ is described in Fig.~5. Analyzing the above density function we get that the
maximal density is achieved at the starting point of the above interval (at the congruent case)
with density $\delta(\mathcal{S}^r(0)) \approx 0.50514$. The density in the endpoint of the above interval is
$\delta(\mathcal{S}^r(h)) \approx 0.23344$.
\begin{figure}[ht]
\centering
\includegraphics[width=8cm]{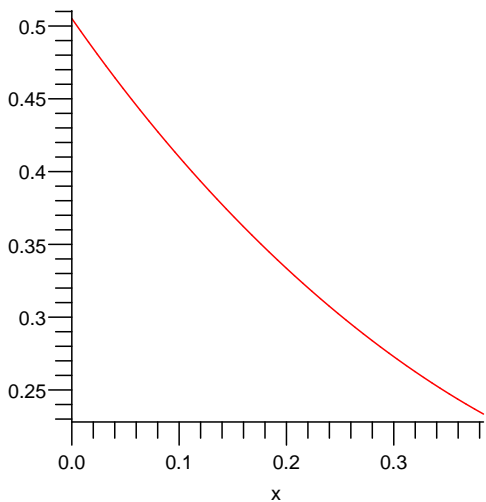}
\caption{Density function of $\delta(\mathcal{S}^r(x))$ where $x \in [0, h]$.}
\label{}
\end{figure}
Finally, we obtain the following
\begin{theorem}
\begin{enumerate}
\item The maximum of the density function $\delta(\mathcal{S}^r(x))$ is attained at the starting
point of the interval $x \in [0,h]$ i.e.
the congruent hyperball packing provides the densest hyperball packing.
Therefore, the density upper bound of the congruent and non-congruent hyperball packings is $\delta(\mathcal{S}^r(0))\approx 0.50514$.
\end{enumerate}
\end{theorem}
The problem of finding the densest hyperball (hypersphere) packing with congruent or non-congruent hyperballs
in $n$-dimensional hyperbolic space ($n\ge3$) is not settled yet. 
At this time the densest hyperball packing with congruent or non-congruent hyperballs is derived by the Coxeter tilings $\{7,3,3\}$ 
(or by the truncated tetrahedron tiling with dihedral angle
$\frac{\pi}{7}$) with density $\approx 0.82251$ but as we have seen, locally there are hyperball packings with larger density than the
B\"or\"oczky-Florian density upper bound for ball and horoball packings (see Theorem 3.2).

We note here, that the discussion of the densest horoball packings in the $n$-dimensional hyperbolic space $n \ge 3$ with horoballs
of different types has not been settled yet as well (see \cite{KSz}, \cite{KSz14}, \cite{Sz12}, \cite{Sz12-2}).

Optimal sphere packings in other homogeneous Thurston geometries represent
another huge class of open mathematical problems. For these non-Euclidean geometries
only very few results are known (e.g. \cite{Sz07-2}, \cite{Sz10}, \cite{Sz13-2}, \cite{Sz14-1}).
Detailed studies are the objective of ongoing research.


\noindent
\footnotesize{Budapest University of Technology and Economics Institute of Mathematics, \\
Department of Geometry, \\
H-1521 Budapest, Hungary. \\
E-mail:~szirmai@math.bme.hu \\
http://www.math.bme.hu/ $^\sim$szirmai}

\end{document}